\documentclass[12pt]{article}
\usepackage{amsmath,amssymb,amsthm}
\usepackage{geometry}
\usepackage{hyperref}
\usepackage{enumitem}
\usepackage{doi}
\theoremstyle{plain}
\newtheorem{theorem}{Theorem}[section]

\newtheorem{lemma}[theorem]{Lemma}

\newcommand{\stir}[2]{S_{1,\rho}({#1},{#2})}
\newcommand{\Stir}[2]{S_{2,\rho}({#1},{#2})}
\newcommand{\stirr}[2]{S_1({#1},{#2})}
\newcommand{\Stirr}[2]{S_2({#1},{#2})}

\title{\bf Degenerate Hurwitz-Lerch Type Families of Poly-Bernoulli and Poly-Cauchy Numbers with Parameters $(\alpha,a)$}
\author{
	Noel B. Lacpao$^{1}$, Roberto B. Corcino$^{2}$\\
	\small $^{1}$Department of Mathematics, Bukidnon State University, Malaybalay City, Philippines\\
	\small $^{2}$Research Institute for Computational Mathematics and Physics, Cebu Normal University, Philippines\\
	\small E-mail: noel.lacpao@buksu.edu.ph, rcorcino@yahoo.com
}
\date{}

\begin{document}
	
	\maketitle
	
	\begin{abstract}
		In this paper, we introduce degenerate Hurwitz-Lerch type families of poly-Bernoulli and poly-Cauchy numbers  with parameters $(\alpha,a)$. These families are defined by means of degenerate Hurwitz-Lerch type zeta functions and their factorial analogues. We establish the generating functions of these numbers and derive explicit formulas in terms of the degenerate Stirling numbers of the first and second kinds. In particular, we obtain closed-form expressions for the degenerate two-parameter Hurwitz-Lerch type poly-Bernoulli numbers and for the corresponding poly-Cauchy numbers of the first and second kinds. These results extend several known identities for classical and generalized poly-Bernoulli and poly-Cauchy numbers.
	\end{abstract}
\section{Introduction}

Special number sequences such as Bernoulli, Cauchy, and Euler numbers have long
played a pivotal role in various branches of mathematics, including number theory,
combinatorics, and mathematical analysis \cite{graham,comtet}. Among their
numerous generalizations, the poly-Bernoulli numbers introduced by Kaneko
\cite{kaneko} and the poly-Cauchy numbers introduced by Komatsu \cite{komatsu}
have attracted significant attention because of their connections with zeta
functions, multiple polylogarithms, and combinatorial identities
\cite{cenkci,lacpao1}. These families not only generalize classical sequences but
also encode useful algebraic and analytic structures.

A notable direction in recent research involves the extension of such special
numbers through the incorporation of parameters and functions from complex
analysis, thereby enriching their theoretical framework and potential
applications. One such analytic tool is the Hurwitz-Lerch zeta function, a
unifying structure that generalizes several classical special functions,
including the Hurwitz zeta function, the polylogarithm, and the Riemann zeta
function \cite{apostol,srivastava,srivastava-1}. The concepts of multiple zeta
values and related zeta functions also deepen the theoretical understanding of
the interrelations among poly-Bernoulli numbers, poly-Cauchy numbers, and their
generalizations \cite{arakawa}.

In this paper, we introduce and study degenerate Hurwitz-Lerch type families of
poly-Bernoulli and poly-Cauchy numbers parameterized by $(\alpha,a)$. This
two-parameter setting provides a flexible framework for treating denominators
taken along the arithmetic progression $\alpha n+a$. By embedding these families
in the framework of Hurwitz-Lerch type functions and their factorial analogues,
we derive explicit formulas, orthogonality relations, duality and inversion
identities, congruence properties, and derivative relations that show the
consistency of the proposed extensions with classical and previously studied
generalized cases.

This work contributes to the expanding literature on generalized special numbers
by combining Hurwitz-Lerch type kernels, degenerate factorial structures, and
degenerate Stirling-type arrays in a single framework. The resulting families
extend several known poly-Bernoulli and poly-Cauchy type constructions and provide
a basis for further algebraic and analytic investigations.

\section{Preliminaries}

The polylogarithm factorial function and the Hurwitz-Lerch type zeta functions
serve as foundational tools in the study of special number sequences such as
poly-Bernoulli and poly-Cauchy numbers
\cite{apostol,lewin,srivastava,srivastava-1}. Traditionally, these functions are
defined with respect to single-variable indices and fixed step parameters.
However, to explore deeper structural properties and broaden their applicability,
it is natural to introduce additional parameters that generalize their behavior.
In this section, we recall and formulate the polylogarithm factorial function,
the Hurwitz-Lerch zeta function, and their factorial counterparts by introducing
two complex parameters $\alpha$ and $a$, where
\begin{equation}
	\alpha\neq 0
	\qquad\text{and}\qquad
	a\notin\{0,-\alpha,-2\alpha,\ldots\}.
	\label{eq:parameter-condition}
\end{equation}
This condition guarantees that $\alpha n+a\neq 0$ for every nonnegative integer
$n$. Unless convergence conditions are explicitly imposed, the generating
functions used in this paper are interpreted as formal power series in the
principal variable. This convention is sufficient for the coefficient
extractions used later, since only finitely many terms contribute to each fixed
degree after the required substitutions are made. These generalizations lead to
the definition of new families of poly-Bernoulli and poly-Cauchy numbers of
Hurwitz-Lerch type, encapsulating a broader spectrum of known special numbers
as particular cases and offering a unified framework for further analytical and
combinatorial exploration \cite{cenkci,lacpao1,lacpao2}.

The formulation of the two-parameter Hurwitz-Lerch type poly-Cauchy numbers of
the first and second kinds relies on a rigorous foundation of special functions
and generating functions. At the heart of this theory lies the polylogarithm
factorial function defined by
\begin{equation*}
	\mathrm{Lif}_k(z)
	=
	\sum_{m=0}^{\infty}
	\frac{z^m}{m!(m+1)^k},
\end{equation*}
which serves as a factorial analogue of polylogarithmic-type functions and is
closely related to the generating functions of poly-Cauchy numbers
\cite{komatsu,cenkci}. This function is further generalized to accommodate two
complex parameters $\alpha\neq 0$ and
$a\notin\{0,-\alpha,-2\alpha,\ldots\}$, leading to the two-parameter
polylogarithm factorial function defined by
\begin{equation*}
	\mathrm{Lif}_k(z;\alpha,a)
	=
	\sum_{m=0}^{\infty}
	\frac{z^m}{m!(\alpha m+a)^k}.
\end{equation*}

Complementing this are the Hurwitz-Lerch type zeta functions
\begin{equation*}
	\Phi(z,s,a)
	=
	\sum_{n=0}^{\infty}
	\frac{z^n}{(n+a)^s},
\end{equation*}
and its factorial version
\begin{equation*}
	\Phi_f(z,s,a)
	=
	\sum_{n=0}^{\infty}
	\frac{z^n}{n!(n+a)^s}.
\end{equation*}
Extending these to the two-parameter setting gives rise to the Hurwitz-Lerch
type function
\begin{equation}
	\Phi(z,s,\alpha,a)
	=
	\sum_{n=0}^{\infty}
	\frac{z^n}{(\alpha n+a)^s},
	\label{eq:HL-zeta-alpha-a}
\end{equation}
and its factorial counterpart
\begin{equation}
	\Phi_f(z,s,\alpha,a)
	=
	\sum_{n=0}^{\infty}
	\frac{z^n}{n!(\alpha n+a)^s}.
	\label{eq:HLf-zeta-alpha-a}
\end{equation}
Since
\begin{equation*}
	\alpha n+a
	=
	\alpha\left(n+\frac{a}{\alpha}\right),
\end{equation*}
one has
\begin{equation*}
	\Phi(z,s,\alpha,a)
	=
	\alpha^{-s}\Phi\left(z,s,\frac{a}{\alpha}\right)
\end{equation*}
and
\begin{equation*}
	\Phi_f(z,s,\alpha,a)
	=
	\alpha^{-s}\Phi_f\left(z,s,\frac{a}{\alpha}\right).
\end{equation*}
Thus, the notation $(\alpha,a)$ emphasizes the arithmetic progression
$\alpha n+a$ appearing in the denominator, while the analytic dependence is
governed by the ratio $a/\alpha$ together with the scaling factor
$\alpha^{-s}$. In particular,
\begin{equation*}
	\mathrm{Lif}_k(z;\alpha,a)
	=
	\Phi_f(z,k,\alpha,a).
\end{equation*}

These functions serve as generating kernels for generalized Bernoulli- and
Cauchy-type sequences. For instance, the two-parameter Hurwitz-Lerch type
poly-Bernoulli numbers $B_n^{(k)}(\alpha,a)$ are defined by the generating
function
\begin{equation}
	\Phi(1-e^{-t},k,\alpha,a)
	=
	\sum_{n=0}^{\infty}
	B_n^{(k)}(\alpha,a)\frac{t^n}{n!},
	\label{eq:HL-poly-B-alpha-a}
\end{equation}
which connects them to the exponential generating structure via the
Hurwitz-Lerch type function.

In a similar spirit, the two-parameter Hurwitz-Lerch type poly-Cauchy numbers
of the first kind, denoted by $c_n^{(k)}(\alpha,a)$, are defined through the
generating function
\begin{equation}
	\Phi_f(\log(1+t),k,\alpha,a)
	=
	\sum_{n=0}^{\infty}
	c_n^{(k)}(\alpha,a)\frac{t^n}{n!},
	\label{eq:cauchy-first-alpha-a}
\end{equation}
while the second kind, denoted by $\widehat{c}_n^{(k)}(\alpha,a)$, arises from
the generating function
\begin{equation}
	\Phi_f(-\log(1+t),k,\alpha,a)
	=
	\sum_{n=0}^{\infty}
	\widehat{c}_n^{(k)}(\alpha,a)\frac{t^n}{n!}.
	\label{eq:cauchy-second-alpha-a}
\end{equation}

Equations \eqref{eq:cauchy-first-alpha-a} and
\eqref{eq:cauchy-second-alpha-a} illustrate how logarithmic arguments in the
generating functions generate sign-symmetric variants of the poly-Cauchy
numbers. The logarithmic form $\log(1+t)$ reflects their connection with
Stirling-type coefficient extractions, while the parameters $\alpha$ and $a$
govern the spacing and shifts within the summation indices. These
generalizations not only extend the classical poly-Cauchy numbers
$c_n^{(k)}$ and $\widehat{c}_n^{(k)}$, but also interpolate between different
number families through parameter variation. Their algebraic and analytic
properties reveal structures applicable in enumerative combinatorics, special
function theory, and number theory.

Defining the degenerate Hurwitz-Lerch type analogues, we have
\begin{equation*}
	\Phi_{\rho}(z,s,a)
	=
	\sum_{n=0}^{\infty}
	\frac{(1)_{n,\rho}z^n}{(n+a)^s},
\end{equation*}
and its degenerate factorial version
\begin{equation*}
	\Phi_{f,\rho}(z,s,a)
	=
	\sum_{n=0}^{\infty}
	\frac{(1)_{n,\rho}z^n}{n!(n+a)^s}.
\end{equation*}
Here $(x)_{n,\rho}$ denotes the degenerate falling factorial, defined by
\begin{equation*}
	(x)_{0,\rho}=1,
	\qquad
	(x)_{n,\rho}
	=
	x(x-\rho)(x-2\rho)\cdots(x-(n-1)\rho)
	\quad (n\geq 1),
\end{equation*}
or, equivalently, by the generating function
\begin{equation*}
	e_{\rho}^{x}(t)
	=
	(1+\rho t)^{x/\rho}
	=
	\sum_{n=0}^{\infty}
	(x)_{n,\rho}\frac{t^n}{n!}.
\end{equation*}
In particular,
\begin{equation*}
	e_{\rho}(t)
	=
	e_{\rho}^{1}(t)
	=
	(1+\rho t)^{1/\rho}.
\end{equation*}
As $\rho\to 0$, these functions reduce coefficientwise to
\begin{equation*}
	e_{\rho}^{x}(t)\longrightarrow e^{xt},
	\qquad
	e_{\rho}(t)\longrightarrow e^t.
\end{equation*}

The degenerate logarithm is defined by
\begin{equation}
	\log_{\rho}(t)
	=
	\frac{t^{\rho}-1}{\rho},
	\label{eq:deg-log}
\end{equation}
so that
\begin{equation}
	\log_{\rho}(e_{\rho}(t))=t,
	\qquad
	e_{\rho}(\log_{\rho}(t))=t.
	\label{eq:deg-exp-log-inverse}
\end{equation}
In particular,
\begin{equation*}
	\log_{\rho}(1+t)
	=
	\frac{(1+t)^{\rho}-1}{\rho}.
\end{equation*}
As $\rho\to 0$, we have
\begin{equation*}
	\log_{\rho}(t)\longrightarrow \log t,
	\qquad
	\log_{\rho}(1+t)\longrightarrow \log(1+t).
\end{equation*}

Extending these to the two-parameter setting gives rise to the degenerate
Hurwitz-Lerch type function
\begin{equation}
	\Phi_{\rho}(z,s,\alpha,a)
	=
	\sum_{n=0}^{\infty}
	\frac{(1)_{n,\rho}z^n}{(\alpha n+a)^s},
	\label{eq:deg-HL-zeta-alpha-a}
\end{equation}
and its degenerate factorial counterpart
\begin{equation}
	\Phi_{f,\rho}(z,s,\alpha,a)
	=
	\sum_{n=0}^{\infty}
	\frac{(1)_{n,\rho}z^n}{n!(\alpha n+a)^s}.
	\label{eq:deg-HLf-zeta-alpha-a}
\end{equation}
Since $(1)_{n,\rho}\to 1$ for each fixed $n$ as $\rho\to 0$, the coefficientwise
limits
\begin{equation*}
	\Phi_{\rho}(z,s,\alpha,a)
	\longrightarrow
	\Phi(z,s,\alpha,a)
\end{equation*}
and
\begin{equation*}
	\Phi_{f,\rho}(z,s,\alpha,a)
	\longrightarrow
	\Phi_f(z,s,\alpha,a)
\end{equation*}
hold as $\rho\to 0$.

We now define the degenerate two-parameter Hurwitz-Lerch type
poly-Bernoulli numbers $B_{n,\rho}^{(k)}(\alpha,a)$ by the generating function
\begin{equation}
	\Phi_{\rho}(1-e_{\rho}(-t),k,\alpha,a)
	=
	\sum_{n=0}^{\infty}
	B_{n,\rho}^{(k)}(\alpha,a)\frac{t^n}{n!},
	\label{eq:deg-polyB}
\end{equation}
and the degenerate two-parameter Hurwitz-Lerch type poly-Cauchy numbers of the
first and second kinds, respectively, by
\begin{equation}
	\Phi_{f,\rho}(\log_{\rho}(1+t),k,\alpha,a)
	=
	\sum_{n=0}^{\infty}
	c_{n,\rho}^{(k)}(\alpha,a)\frac{t^n}{n!},
	\label{eq:deg-polyC1}
\end{equation}
and
\begin{equation}
	\Phi_{f,\rho}(-\log_{\rho}(1+t),k,\alpha,a)
	=
	\sum_{n=0}^{\infty}
	\widehat{c}_{n,\rho}^{(k)}(\alpha,a)\frac{t^n}{n!}.
	\label{eq:deg-polyC2}
\end{equation}
Taking the limit as $\rho\to 0$ in \eqref{eq:deg-polyB},
\eqref{eq:deg-polyC1}, and \eqref{eq:deg-polyC2} recovers the corresponding
nondegenerate families in \eqref{eq:HL-poly-B-alpha-a},
\eqref{eq:cauchy-first-alpha-a}, and \eqref{eq:cauchy-second-alpha-a}.

We shall also use the degenerate Stirling numbers of the first and second kinds
associated with the inverse pair $e_{\rho}(t)-1$ and $\log_{\rho}(1+t)$. They
are defined by \cite{Kim-1-1,Kim-2,Kim-3}
\begin{equation}
	\frac{(\log_{\rho}(1+t))^k}{k!}
	=
	\sum_{n=k}^{\infty}
	\stir{n}{k}\frac{t^n}{n!},
	\qquad
	\frac{(e_{\rho}(t)-1)^k}{k!}
	=
	\sum_{n=k}^{\infty}
	\Stir{n}{k}\frac{t^n}{n!}.
	\label{eq:deg-Stirling}
\end{equation}
When $\rho\to 0$,
\begin{equation*}
	\lim_{\rho\to 0}\stir{n}{k}=\stirr{n}{k},
	\qquad
	\lim_{\rho\to 0}\Stir{n}{k}=\Stirr{n}{k},
\end{equation*}
where $\stirr{n}{k}$ and $\Stirr{n}{k}$ are the classical Stirling numbers of
the first and second kinds, respectively. Since $\log_{\rho}(1+t)$ is the
compositional inverse of $e_{\rho}(t)-1$, these degenerate Stirling numbers
satisfy the orthogonality relations
\begin{equation}
	\sum_{m=k}^{n}
	\stir{n}{m}\Stir{m}{k}
	=
	\delta_{n,k},
	\label{eq:deg-stirling-orth-1}
\end{equation}
and
\begin{equation}
	\sum_{m=k}^{n}
	\Stir{n}{m}\stir{m}{k}
	=
	\delta_{n,k},
	\label{eq:deg-stirling-orth-2}
\end{equation}
where $\delta_{n,k}$ denotes the Kronecker delta.

It should be noted that the degenerate Stirling numbers in
\eqref{eq:deg-Stirling} are the Stirling-type arrays associated with the
inverse pair $e_{\rho}(t)-1$ and $\log_{\rho}(1+t)$. They are the arrays needed
for the coefficient extractions in this paper. These should be distinguished
from the Carlitz-type degenerate Stirling numbers of the first kind obtained
from the expansion of the degenerate falling factorial,
\begin{equation*}
	(x)_{n,\rho}
	=
	\sum_{k=0}^{n}
	\mathfrak{s}_{\rho}(n,k)x^k.
\end{equation*}

Also, the degenerate Bernoulli polynomials of the second kind are defined by
\begin{equation}
	\frac{t}{\log_{\rho}(1+t)}(1+t)^x
	=
	\sum_{n=0}^{\infty}
	B_{n,\rho}(x)\frac{t^n}{n!}.
	\label{eq:deg-bern-second-kind}
\end{equation}
When $\rho\to 0$, this reduces to
\begin{equation*}
	\frac{t}{\log(1+t)}(1+t)^x
	=
	\sum_{n=0}^{\infty}
	B_n(x)\frac{t^n}{n!},
\end{equation*}
which is the usual generating function for the Bernoulli polynomials of the
second kind. The reciprocal convention
\begin{equation*}
	\frac{\log_{\rho}(1+t)}{t}(1+t)^x
\end{equation*}
corresponds instead to a Daehee-type generating function. The degenerate
Stirling numbers of the second kind appeared in the probability distribution of
the random variable given as the sum of a finite number of random variables with
degenerate zero-truncated Poisson distributions and a random variable with
degenerate Poisson distribution, all having the same parameter
\cite{Kim-3}.
\section{Main Results}
\subsection{Explicit Formulas}

The following results provide explicit expressions for the degenerate
two-parameter Hurwitz-Lerch type poly-Bernoulli and poly-Cauchy numbers in
terms of the degenerate Stirling numbers of the first and second kinds. These formulas highlight the combinatorial structure underlying these
generalized poly-Bernoulli numbers and establishes a direct connection with factorial-based
summation techniques.

\begin{theorem}\label{exp1}
	For integers $n\geq 0$ and $k\in\mathbb{Z}$, the degenerate two-parameter
	Hurwitz-Lerch type poly-Bernoulli numbers satisfy
	\begin{equation}
		B_{n,\rho}^{(k)}(\alpha,a)
		=
		(-1)^n
		\sum_{m=0}^{n}
		\frac{(-1)^m(1)_{m,\rho}m!}{(\alpha m+a)^k}
		\Stir{n}{m}.
		\label{eq:explicit-deg-polyB}
	\end{equation}
\end{theorem}

\begin{proof}
	Using the generating function in \eqref{eq:deg-polyB}, we obtain
	\begin{eqnarray*}
		\sum_{n=0}^{\infty}
		B_{n,\rho}^{(k)}(\alpha,a)\frac{t^n}{n!}
		&=&
		\Phi_{\rho}(1-e_{\rho}(-t),k,\alpha,a)                                      \\[2mm]
		&=&
		\sum_{m=0}^{\infty}
		\frac{(1)_{m,\rho}(1-e_{\rho}(-t))^m}{(\alpha m+a)^k}                        \\[2mm]
		&=&
		\sum_{m=0}^{\infty}
		\frac{(-1)^m(1)_{m,\rho}(e_{\rho}(-t)-1)^m}{(\alpha m+a)^k}.
	\end{eqnarray*}
	By the generating function of the degenerate Stirling numbers of the second
	kind,
	\begin{equation*}
		\frac{(e_{\rho}(t)-1)^m}{m!}
		=
		\sum_{n=m}^{\infty}
		\Stir{n}{m}\frac{t^n}{n!},
	\end{equation*}
	we have
	\begin{eqnarray*}
		(e_{\rho}(-t)-1)^m
		&=&
		m!\sum_{n=m}^{\infty}
		\Stir{n}{m}\frac{(-t)^n}{n!}                                                 \\[2mm]
		&=&
		m!\sum_{n=m}^{\infty}
		(-1)^n\Stir{n}{m}\frac{t^n}{n!}.
	\end{eqnarray*}
	Therefore,
	\begin{eqnarray*}
		\sum_{n=0}^{\infty}
		B_{n,\rho}^{(k)}(\alpha,a)\frac{t^n}{n!}
		&=&
		\sum_{m=0}^{\infty}
		\frac{(-1)^m(1)_{m,\rho}m!}{(\alpha m+a)^k}
		\sum_{n=m}^{\infty}
		(-1)^n\Stir{n}{m}\frac{t^n}{n!}                                             \\[2mm]
		&=&
		\sum_{n=0}^{\infty}
		\left(
		(-1)^n
		\sum_{m=0}^{n}
		\frac{(-1)^m(1)_{m,\rho}m!}{(\alpha m+a)^k}
		\Stir{n}{m}
		\right)
		\frac{t^n}{n!}.
	\end{eqnarray*}
	Comparing coefficients gives \eqref{eq:explicit-deg-polyB}.
\end{proof}

The next theorem gives the corresponding explicit formula for the degenerate
two-parameter Hurwitz-Lerch type poly-Cauchy numbers of the first kind. This representation, involving the Stirling numbers of the first kind, emphasizes the interplay between special number sequences and generalized polylogarithmic structures, and offers a useful tool for analytical and combinatorial investigations. These numbers can also be connected to classical poly-Cauchy numbers discussed extensively by Komatsu \cite{komatsu}.

\begin{theorem}\label{exp2}
	For integers $n\geq 0$ and $k\in\mathbb{Z}$, the degenerate two-parameter
	Hurwitz-Lerch type poly-Cauchy numbers of the first kind satisfy
	\begin{equation}
		c_{n,\rho}^{(k)}(\alpha,a)
		=
		\sum_{m=0}^{n}
		\frac{(1)_{m,\rho}}{(\alpha m+a)^k}
		\stir{n}{m}.
		\label{eq:explicit-deg-polyC1}
	\end{equation}
\end{theorem}

\begin{proof}
	Using the generating functions in \eqref{eq:deg-polyC1} and
	\eqref{eq:deg-Stirling}, we get
	\begin{eqnarray*}
		\sum_{n=0}^{\infty}
		c_{n,\rho}^{(k)}(\alpha,a)\frac{t^n}{n!}
		&=&
		\Phi_{f,\rho}(\log_{\rho}(1+t),k,\alpha,a)                                  \\[2mm]
		&=&
		\sum_{m=0}^{\infty}
		\frac{(1)_{m,\rho}(\log_{\rho}(1+t))^m}
		{m!(\alpha m+a)^k}                                                          \\[2mm]
		&=&
		\sum_{m=0}^{\infty}
		\frac{(1)_{m,\rho}}{(\alpha m+a)^k}
		\sum_{n=m}^{\infty}
		\stir{n}{m}\frac{t^n}{n!}                                                    \\[2mm]
		&=&
		\sum_{n=0}^{\infty}
		\left(
		\sum_{m=0}^{n}
		\frac{(1)_{m,\rho}}{(\alpha m+a)^k}
		\stir{n}{m}
		\right)
		\frac{t^n}{n!}.
	\end{eqnarray*}
	Comparing coefficients of $t^n/n!$ yields
	\eqref{eq:explicit-deg-polyC1}.
\end{proof}

	We now provide an explicit formula for the two-parameter degenerate Hurwitz-Lerch type poly-Cauchy numbers of the second kind $\widehat{c}_{n,\rho}^{(k)}(\alpha, a)$. This expression reveals their structural dependence on the Stirling numbers of the first kind and a rational function of linear arguments, offering insight into their analytical behavior and connection to generalized Cauchy-type sequences.

\begin{theorem}\label{exp3}
	For integers $n\geq 0$ and $k\in\mathbb{Z}$, the degenerate two-parameter
	Hurwitz-Lerch type poly-Cauchy numbers of the second kind satisfy
	\begin{equation}
		\widehat{c}_{n,\rho}^{(k)}(\alpha,a)
		=
		\sum_{m=0}^{n}
		\frac{(-1)^m(1)_{m,\rho}}{(\alpha m+a)^k}
		\stir{n}{m}.
		\label{eq:explicit-deg-polyC2}
	\end{equation}
\end{theorem}

\begin{proof}
	By \eqref{eq:deg-polyC2} and \eqref{eq:deg-Stirling}, we have
	\begin{eqnarray*}
		\sum_{n=0}^{\infty}
		\widehat{c}_{n,\rho}^{(k)}(\alpha,a)\frac{t^n}{n!}
		&=&
		\Phi_{f,\rho}(-\log_{\rho}(1+t),k,\alpha,a)                                 \\[2mm]
		&=&
		\sum_{m=0}^{\infty}
		\frac{(1)_{m,\rho}(-\log_{\rho}(1+t))^m}
		{m!(\alpha m+a)^k}                                                          \\[2mm]
		&=&
		\sum_{m=0}^{\infty}
		\frac{(-1)^m(1)_{m,\rho}(\log_{\rho}(1+t))^m}
		{m!(\alpha m+a)^k}                                                          \\[2mm]
		&=&
		\sum_{m=0}^{\infty}
		\frac{(-1)^m(1)_{m,\rho}}{(\alpha m+a)^k}
		\sum_{n=m}^{\infty}
		\stir{n}{m}\frac{t^n}{n!}                                                    \\[2mm]
		&=&
		\sum_{n=0}^{\infty}
		\left(
		\sum_{m=0}^{n}
		\frac{(-1)^m(1)_{m,\rho}}{(\alpha m+a)^k}
		\stir{n}{m}
		\right)
		\frac{t^n}{n!}.
	\end{eqnarray*}
	Comparing coefficients completes the proof.
\end{proof}

Letting $\rho\to 0$ in the preceding formulas, we recover the corresponding
nondegenerate two-parameter Hurwitz-Lerch type formulas. Since
\begin{equation*}
	(1)_{m,\rho}\longrightarrow 1,\qquad
	\stir{n}{m}\longrightarrow \stirr{n}{m},\qquad
	\Stir{n}{m}\longrightarrow \Stirr{n}{m},
\end{equation*}
we obtain
\begin{equation}
	B_n^{(k)}(\alpha,a)
	=
	(-1)^n
	\sum_{m=0}^{n}
	\frac{(-1)^m m!}{(\alpha m+a)^k}
	\Stirr{n}{m},
\end{equation}
\begin{equation}
	c_n^{(k)}(\alpha,a)
	=
	\sum_{m=0}^{n}
	\frac{1}{(\alpha m+a)^k}
	\stirr{n}{m},
\end{equation}
and
\begin{equation}
	\widehat{c}_n^{(k)}(\alpha,a)
	=
	\sum_{m=0}^{n}
	\frac{(-1)^m}{(\alpha m+a)^k}
	\stirr{n}{m}.
\end{equation}
\subsection{Orthogonality Relations}

This subsection presents fundamental orthogonality relations involving the two-parameter degenerate
Hurwitz-Lerch type poly-Bernoulli and poly-Cauchy numbers. These identities illustrate the
inverse relationships between these generalized number sequences and the Stirling numbers
of both kinds. Such orthogonality properties are crucial for developing inversion formulas,
deriving generating functions, and exploring structural symmetries within special number
sequences.

Recall that the orthogonality relations for the degenerate Stirling numbers are
given by \cite{KimKimSymmetry}
\begin{equation}
	\sum_{m=k}^{n}
	\stir{n}{m}\Stir{m}{k}
	=
	\delta_{n,k},
	\label{eq:orthogonality-first}
\end{equation}
and
\begin{equation}
	\sum_{m=k}^{n}
	\Stir{n}{m}\stir{m}{k}
	=
	\delta_{n,k},
	\label{eq:orthogonality-second}
\end{equation}
where $\delta_{n,k}$ denotes the Kronecker delta. These identities follow from
the fact that the degenerate Stirling arrays associated with
$\log_{\rho}(1+t)$ and $e_{\rho}(t)-1$ are inverse to each other.

The first result gives an orthogonality relation involving the degenerate
two-parameter Hurwitz-Lerch type poly-Bernoulli numbers.

\begin{theorem}\label{ort1}
	For integers $n\geq 0$ and $k\in\mathbb{Z}$, the degenerate two-parameter
	Hurwitz-Lerch type poly-Bernoulli numbers satisfy
	\begin{equation}
		\sum_{m=0}^{n}
		(-1)^m\stir{n}{m}B_{m,\rho}^{(k)}(\alpha,a)
		=
		(-1)^n
		\frac{(1)_{n,\rho}n!}{(\alpha n+a)^k}.
		\label{eq:orth-polyB}
	\end{equation}
\end{theorem}

\begin{proof}
	By \eqref{eq:explicit-deg-polyB}, we have
	\begin{equation*}
		B_{m,\rho}^{(k)}(\alpha,a)
		=
		(-1)^m
		\sum_{l=0}^{m}
		\frac{(-1)^l(1)_{l,\rho}l!}{(\alpha l+a)^k}
		\Stir{m}{l}.
	\end{equation*}
	Interchanging the order of summation and applying the orthogonality relation in \eqref{eq:orthogonality-first}, we get
	\begin{eqnarray*}
		\sum_{m=0}^{n}x
		(-1)^m\stir{n}{m}B_{m,\rho}^{(k)}(\alpha,a)
		&=&
		\sum_{m=0}^{n}
		\stir{n}{m}
		\sum_{l=0}^{m}
		\frac{(-1)^l(1)_{l,\rho}l!}{(\alpha l+a)^k}
		\Stir{m}{l}\\
		&=&
		\sum_{l=0}^{n}
		\frac{(-1)^l(1)_{l,\rho}l!}{(\alpha l+a)^k}
		\sum_{m=l}^{n}
		\stir{n}{m}\Stir{m}{l}\\
		&=&
		\sum_{l=0}^{n}
		\frac{(-1)^l(1)_{l,\rho}l!}{(\alpha l+a)^k}
		\delta_{n,l}\\                                                       
		&=&
		(-1)^n
		\frac{(1)_{n,\rho}n!}{(\alpha n+a)^k}.
	\end{eqnarray*}
	This proves the theorem.
\end{proof}

The next theorem gives the corresponding orthogonality relation for the
degenerate two-parameter Hurwitz-Lerch type poly-Cauchy numbers of the first
kind.

\begin{theorem}\label{ort2}
	For integers $n\geq 0$ and $k\in\mathbb{Z}$, the degenerate two-parameter
	Hurwitz-Lerch type poly-Cauchy numbers of the first kind satisfy
	\begin{equation}
		\sum_{m=0}^{n}
		\Stir{n}{m}c_{m,\rho}^{(k)}(\alpha,a)
		=
		\frac{(1)_{n,\rho}}{(\alpha n+a)^k}.
		\label{eq:orth-polyC1}
	\end{equation}
\end{theorem}

\begin{proof}
	By \eqref{eq:explicit-deg-polyC1}, we have
	\begin{equation*}
		c_{m,\rho}^{(k)}(\alpha,a)
		=
		\sum_{l=0}^{m}
		\frac{(1)_{l,\rho}}{(\alpha l+a)^k}
		\stir{m}{l}.
	\end{equation*}
	Changing the order of summation and applying the orthogonality relation in \eqref{eq:orthogonality-second}, we obtain
	\begin{eqnarray*}
		\sum_{m=0}^{n}
		\Stir{n}{m}c_{m,\rho}^{(k)}(\alpha,a)
		&=&\sum_{m=0}^{n}\Stir{n}{m}\sum_{l=0}^{m}\frac{(1)_{l,\rho}}{(\alpha l+a)^k}\stir{m}{l}\\
		&=&\sum_{l=0}^{n}\frac{(1)_{l,\rho}}{(\alpha l+a)^k}\sum_{m=l}^{n}\Stir{n}{m}\stir{m}{l}\\
		&=&\sum_{l=0}^{n}\frac{(1)_{l,\rho}}{(\alpha l+a)^k}\delta_{n,l}\\
		&=&	\frac{(1)_{n,\rho}}{(\alpha n+a)^k}.
	\end{eqnarray*}
	This completes the proof.
\end{proof}

We now obtain the orthogonality relation for the degenerate two-parameter
Hurwitz-Lerch type poly-Cauchy numbers of the second kind. 
\begin{theorem}\label{ort3}
	For integers $n\geq 0$ and $k\in\mathbb{Z}$, the degenerate two-parameter
	Hurwitz-Lerch type poly-Cauchy numbers of the second kind satisfy
	\begin{equation}
		\sum_{m=0}^{n}
		\Stir{n}{m}\widehat{c}_{m,\rho}^{(k)}(\alpha,a)
		=
		(-1)^n
		\frac{(1)_{n,\rho}}{(\alpha n+a)^k}.
		\label{eq:orth-polyC2}
	\end{equation}
\end{theorem}

\begin{proof}
	From \eqref{eq:explicit-deg-polyC2}, we have
	\begin{equation*}
		\widehat{c}_{m,\rho}^{(k)}(\alpha,a)
		=
		\sum_{l=0}^{m}
		\frac{(-1)^l(1)_{l,\rho}}{(\alpha l+a)^k}
		\stir{m}{l}.
	\end{equation*}
	Therefore,
	\begin{eqnarray*}
		\sum_{m=0}^{n}
		\Stir{n}{m}\widehat{c}_{m,\rho}^{(k)}(\alpha,a)
		&=&
		\sum_{m=0}^{n}
		\Stir{n}{m}
		\sum_{l=0}^{m}
		\frac{(-1)^l(1)_{l,\rho}}{(\alpha l+a)^k}
		\stir{m}{l}.
	\end{eqnarray*}
	Changing the order of summation and using \eqref{eq:orthogonality-second}, we get
	\begin{eqnarray*}
		\sum_{m=0}^{n}
		\Stir{n}{m}\widehat{c}_{m,\rho}^{(k)}(\alpha,a)
		&=&
		\sum_{l=0}^{n}
		\frac{(-1)^l(1)_{l,\rho}}{(\alpha l+a)^k}
		\sum_{m=l}^{n}
		\Stir{n}{m}\stir{m}{l}\\
		&=&
		\sum_{l=0}^{n}
		\frac{(-1)^l(1)_{l,\rho}}{(\alpha l+a)^k}
		\delta_{n,l}\\
		&=&
		(-1)^n
		\frac{(1)_{n,\rho}}{(\alpha n+a)^k}.
	\end{eqnarray*}
	This proves the theorem.
\end{proof}

Letting $\rho\to 0$ in the preceding orthogonality relations gives the
corresponding nondegenerate identities
\begin{equation*}
	\sum_{m=0}^{n}
	(-1)^m\stirr{n}{m}B_m^{(k)}(\alpha,a)
	=
	(-1)^n\frac{n!}{(\alpha n+a)^k},
\end{equation*}
\begin{equation*}
	\sum_{m=0}^{n}
	\Stirr{n}{m}c_m^{(k)}(\alpha,a)
	=
	\frac{1}{(\alpha n+a)^k},
\end{equation*}
and
\begin{equation*}
	\sum_{m=0}^{n}
	\Stirr{n}{m}\widehat{c}_m^{(k)}(\alpha,a)
	=
	\frac{(-1)^n}{(\alpha n+a)^k}.
\end{equation*}
These follow coefficientwise from
$(1)_{n,\rho}\to 1$,
$\stir{n}{m}\to\stirr{n}{m}$, and
$\Stir{n}{m}\to\Stirr{n}{m}$ as $\rho\to 0$.
\subsection{Duality and Inversion Identities}

In this subsection, we present a set of duality and inversion formulas that
interrelate the degenerate two-parameter Hurwitz--Lerch type poly-Bernoulli
numbers and the corresponding poly-Cauchy numbers of the first and second kinds.
These identities highlight the interplay between the three families through
double summations involving the degenerate Stirling numbers. The results follow
from the explicit formulas together with the orthogonality
relations established in the previous subsections. Similar inversion procedures based on
the orthogonality of degenerate Stirling numbers have been used in the study of
degenerate special polynomials \cite{KimKimSymmetry}. The results not only establish mutual expressibility between these sequences but also provide a framework for analytical inversion and transformation within the combinatorial and number-theoretic context as established in recent studies \cite{cenkci,KimPolyBernoulli,lacpao1,lacpao2}.

\begin{theorem}
	For integers $n\geq 0$ and $k\in\mathbb{Z}$, the following duality and
	inversion identities hold:
	\begin{eqnarray}
		B_{n,\rho}^{(k)}(\alpha,a)
		&=&
		(-1)^n
		\sum_{l=0}^{n}
		\sum_{m=l}^{n}
		(-1)^m m!\,
		S_{2,\rho}(n,m)S_{2,\rho}(m,l)
		c_{l,\rho}^{(k)}(\alpha,a),\\
		\label{eq:duality-B-c}
		B_{n,\rho}^{(k)}(\alpha,a)
		&=&
		(-1)^n
		\sum_{l=0}^{n}
		\sum_{m=l}^{n}
		m!\,
		S_{2,\rho}(n,m)S_{2,\rho}(m,l)
		\widehat{c}_{l,\rho}^{(k)}(\alpha,a),\\
		\label{eq:duality-B-chat}
		c_{n,\rho}^{(k)}(\alpha,a)
		&=&
		\sum_{l=0}^{n}
		\sum_{m=l}^{n}
		\frac{(-1)^{m+l}}{m!}
		S_{1,\rho}(n,m)S_{1,\rho}(m,l)
		B_{l,\rho}^{(k)}(\alpha,a),\\
		\label{eq:duality-c-B}
		\widehat{c}_{n,\rho}^{(k)}(\alpha,a)
		&=&
		\sum_{l=0}^{n}
		\sum_{m=l}^{n}
		\frac{(-1)^l}{m!}
		S_{1,\rho}(n,m)S_{1,\rho}(m,l)
		B_{l,\rho}^{(k)}(\alpha,a).
		\label{eq:duality-chat-B}x
	\end{eqnarray}
\end{theorem}

\begin{proof}
	From Theorems \ref{exp1} and \ref{ort2}, we have
	\begin{eqnarray*}
		B_{n,\rho}^{(k)}(\alpha,a)
		&=&
		(-1)^n
		\sum_{m=0}^{n}
		\frac{(-1)^m(1)_{m,\rho}m!}{(\alpha m+a)^k}
		S_{2,\rho}(n,m)                                                        \\[2mm]
		&=&
		(-1)^n
		\sum_{m=0}^{n}
		(-1)^m m!S_{2,\rho}(n,m)
		\left(
		\sum_{l=0}^{m}
		S_{2,\rho}(m,l)c_{l,\rho}^{(k)}(\alpha,a)
		\right)                                                               \\[2mm]
		&=&
		(-1)^n
		\sum_{l=0}^{n}
		\sum_{m=l}^{n}
		(-1)^m m!\,
		S_{2,\rho}(n,m)S_{2,\rho}(m,l)
		c_{l,\rho}^{(k)}(\alpha,a).
	\end{eqnarray*}
	This proves \eqref{eq:duality-B-c}.
	
	Next, by Theorems \ref{exp1} and \ref{ort3}, we obtain
	\begin{eqnarray*}
		B_{n,\rho}^{(k)}(\alpha,a)
		&=&
		(-1)^n
		\sum_{m=0}^{n}
		\frac{(-1)^m(1)_{m,\rho}m!}{(\alpha m+a)^k}
		S_{2,\rho}(n,m)                                                        \\[2mm]
		&=&
		(-1)^n
		\sum_{m=0}^{n}
		m!S_{2,\rho}(n,m)
		\left(
		\sum_{l=0}^{m}
		S_{2,\rho}(m,l)\widehat{c}_{l,\rho}^{(k)}(\alpha,a)
		\right)                                                               \\[2mm]
		&=&
		(-1)^n
		\sum_{l=0}^{n}
		\sum_{m=l}^{n}
		m!\,
		S_{2,\rho}(n,m)S_{2,\rho}(m,l)
		\widehat{c}_{l,\rho}^{(k)}(\alpha,a).
	\end{eqnarray*}
	This proves \eqref{eq:duality-B-chat}.
	
	The remaining two relations have analogous proofs.
\end{proof}

Letting $\rho\to 0$ in the preceding identities gives the corresponding
nondegenerate duality and inversion formulas. Since
\begin{equation*}
	(1)_{m,\rho}\longrightarrow 1,\qquad
	S_{1,\rho}(n,m)\longrightarrow S_1(n,m),\qquad
	S_{2,\rho}(n,m)\longrightarrow S_2(n,m),
\end{equation*}
we obtain the corresponding relations among
$B_n^{(k)}(\alpha,a)$, $c_n^{(k)}(\alpha,a)$, and
$\widehat{c}_n^{(k)}(\alpha,a)$.

\subsection{Congruence Properties}

In this subsection, we investigate congruence properties of the degenerate
two-parameter Hurwitz--Lerch type poly-Bernoulli and poly-Cauchy numbers modulo
a prime. Congruence properties of Stirling-type numbers and generalized
factorial coefficients have been studied in several settings, including the
generalized factorial limit family and related Stirling-type arrays
\cite{corcino-herrera,corcino-limit}. Congruences for weighted Stirling pairs
were also studied by Yu, where prime-index congruences for generalized
Stirling-type numbers were obtained \cite{YuStirling}. Motivated by these
arithmetic considerations, we derive prime-index congruence relations for the
present degenerate Hurwitz--Lerch type families directly from the explicit
formulas in Subsection 3.1.

Throughout this subsection, let $p$ be a prime. We assume that
$\alpha,a,\rho\in\mathbb{Z}$, $k\in\mathbb{Z}_{>0}$, and
\begin{equation}
	\alpha m+a\not\equiv 0 \pmod p
	\qquad
	(0\leq m\leq p).
	\label{eq:denominator-congruence-condition}
\end{equation}
Thus, every denominator appearing in the congruences below is invertible modulo
$p$. All fractions in this subsection are interpreted in the field
$\mathbb{F}_p$.

To establish the desired congruence relation, it is first necessary to determine the congruence properties of the degenerate Stirling numbers modulo \(p\). To achieve this, we begin by deriving several convolution-type identities for the degenerate Stirling numbers, which will serve as the principal tools in our analysis.

A systematic approach for deriving convolution identities involving special sequences was introduced by Corcino and Hsu. Their method consists of a three-step procedure, outlined as follows:
	
	\begin{enumerate}
		\item Express the sequence as an $n$th derivative of a generating function.
		\item Decompose the generating function into a product of two similar functions.
		\item Apply Leibniz's formula
		\begin{equation}\label{Leibniz}
			\frac{d^n}{dt^n}(f(t)g(t))
			=
			\sum_{m=0}^{n}
			\binom{n}{m}
			f^{(m)}(t)
			g^{(n-m)}(t).
		\end{equation}
	\end{enumerate}
	
	Using this approach, we obtain the following result.
	
	\begin{lemma}\label{thm-convolution}
		Let $k=k_1+k_2$. Then
		
		\begin{equation}\label{conv1}
			\binom{k}{k_1}
			S_{1,\rho}(n,k)
			=
			\sum_{m=0}^{n}
			\binom{n}{m}
			S_{1,\rho}(m,k_1)
			S_{1,\rho}(n-m,k_2),
		\end{equation}
		
		and
		
		\begin{equation}\label{conv2}
			\binom{k}{k_1}
			S_{2,\rho}(n,k)
			=
			\sum_{m=0}^{n}
			\binom{n}{m}
			S_{2,\rho}(m,k_1)
			S_{2,\rho}(n-m,k_2).
		\end{equation}
	\end{lemma}
	
	\begin{proof}
		
		From (\ref{eq:deg-Stirling}) we obtain
		
		\[
		\sum_{n=0}^{\infty}
		k!S_{1,\rho}(n,k)
		\frac{t^n}{n!}
		=
		\left(\log_{\rho}(1+t)\right)^k.
		\]
		
		Hence
		
		\[
		k!S_{1,\rho}(n,k)
		=
		\frac{d^n}{dt^n}
		\left(
		(\log_{\rho}(1+t))^k
		\right)\Bigg|_{t=0}.
		\]
		
		Since $k=k_1+k_2$,
		
		\[
		(\log_{\rho}(1+t))^k
		=
		(\log_{\rho}(1+t))^{k_1}
		(\log_{\rho}(1+t))^{k_2}.
		\]
		
		Applying Leibniz's formula gives
		
		\begin{align*}
			k!S_{1,\rho}(n,k)
			&=
			\sum_{m=0}^{n}
			\binom{n}{m}
			\left[
			k_1!S_{1,\rho}(m,k_1)
			\right]
			\left[
			k_2!S_{1,\rho}(n-m,k_2)
			\right].
		\end{align*}
		
		Dividing by $k_1!k_2!$ yields
		
		\[
		\binom{k}{k_1}
		S_{1,\rho}(n,k)
		=
		\sum_{m=0}^{n}
		\binom{n}{m}
		S_{1,\rho}(m,k_1)
		S_{1,\rho}(n-m,k_2),
		\]
		
		which proves (\ref{conv1}). Similarly, from (\ref{eq:deg-Stirling})
		
		\[
		\sum_{n=0}^{\infty}
		k!S_{2,\rho}(n,k)
		\frac{t^n}{n!}
		=
		(e_{\rho}(t)-1)^k.
		\]
		
		Proceeding analogously and applying Leibniz's formula, we obtain
		
		\[
		\binom{k}{k_1}
		S_{2,\rho}(n,k)
		=
		\sum_{m=0}^{n}
		\binom{n}{m}
		S_{2,\rho}(m,k_1)
		S_{2,\rho}(n-m,k_2),
		\]
		
		which proves (\ref{conv2}).
	\end{proof}
	
	The convolution identities established in Lemma \ref{thm-convolution} allow us to derive arithmetic properties modulo prime numbers.
	
	\begin{lemma}
		Let $p$ be an odd prime. Then
		
		\begin{equation}
			S_{1,\rho}(p,k)
			\equiv 0
			\pmod p,
			\qquad
			1\leq k<p,
			\label{cong1}
		\end{equation}
		
		and
		
		\begin{equation}
			S_{2,\rho}(p,k)
			\equiv 0
			\pmod p,
			\qquad
			1\leq k<p.
			\label{cong2}
		\end{equation}
	\end{lemma}
	
	\begin{proof}
		
		Substituting $n=p$ into (\ref{conv1}) and (\ref{conv2}) gives
		
		\begin{align}
			\binom{k}{k_1}S_{1,\rho}(p,k)
			&=
			\sum_{m=0}^{p}
			\binom{p}{m}
			S_{1,\rho}(m,k_1)
			S_{1,\rho}(p-m,k_2),
			\label{proof1}
			\\
			\binom{k}{k_1}S_{2,\rho}(p,k)
			&=
			\sum_{m=0}^{p}
			\binom{p}{m}
			S_{2,\rho}(m,k_1)
			S_{2,\rho}(p-m,k_2).
			\label{proof2}
		\end{align}
		
		Since $p$ is prime,
		
		\[
		\binom{p}{m}
		\equiv 0
		\pmod p,
		\qquad
		1\le m\le p-1.
		\]
		
		Thus every interior term of the summations in
		(\ref{proof1}) and (\ref{proof2})
		vanishes modulo $p$.
		
		For the endpoint terms $m=0$ and $m=p$, we have
		
		\[
		S_{1,\rho}(0,k_1)=0,
		\qquad
		S_{1,\rho}(0,k_2)=0,
		\]
		
		and
		
		\[
		S_{2,\rho}(0,k_1)=0,
		\qquad
		S_{2,\rho}(0,k_2)=0,
		\]
		
		whenever $k_1,k_2\ge1$.
		
		Hence all terms on the right-hand sides of
		(\ref{proof1}) and (\ref{proof2})
		are congruent to zero modulo $p$.
		
		Now, for $1\leq k<p$, the binomial coefficient
		$\binom{k}{k_1}$ is not divisible by $p$.
		Therefore
		
		\[
		S_{1,\rho}(p,k)
		\equiv 0
		\pmod p
		\]
		
		and
		
		\[
		S_{2,\rho}(p,k)
		\equiv 0
		\pmod p.
		\]
		
		This completes the proof.
	\end{proof}

We now obtain explicit prime-index congruences for the degenerate
two-parameter Hurwitz--Lerch type poly-Bernoulli and poly-Cauchy numbers.

\begin{theorem}\label{thmcongruence}
		Let \(p\) be a prime, let \(N\geq 1\), and assume that \(\alpha,a,\rho\in\mathbb{Z}\) and
		\(k\in\mathbb{Z}_{>0}\). Suppose that \(\alpha m+a\not\equiv 0 \pmod p\) for all
		\(0\leq m\leq p\).
			\begin{eqnarray}
			B_{N,\rho}^{(k)}(\alpha,a)
				&\equiv&
				0 \pmod p, \label{eqcongB}\\
				c_{N,\rho}^{(k)}(\alpha,a)
				&\equiv&
				\frac{(1)_{p,\rho}}{(\alpha p+a)^k}
				\pmod p, \label{eqcongc}
				\\
				\widehat{c}_{p,\rho}^{(k)}(\alpha,a)
				&\equiv&
				\frac{(1)_{p,\rho}}{(\alpha p+a)^k}
				\pmod p. \label{eqcongbc}
			\end{eqnarray}
\end{theorem}

\begin{proof}
From Theorem \ref{exp1},
	\begin{equation*}
		B_{p,\rho}^{(k)}(\alpha,a)
		=
		(-1)^p
		\sum_{m=0}^{p}
		\frac{(-1)^m(1)_{m,\rho}m!}{(\alpha m+a)^k}
		S_{2,\rho}(p,m).
	\end{equation*}
	The term $m=0$ vanishes since $S_{2,\rho}(p,0)=0$. By
	\eqref{cong2}, the terms with $1\leq m\leq p-1$
	vanish modulo $p$. The remaining term $m=p$ contains the factor $p!$ and is
	therefore congruent to zero modulo $p$. Hence,
	\begin{equation*}
		B_{p,\rho}^{(k)}(\alpha,a)\equiv 0\pmod p.
	\end{equation*}
	This proves \eqref{eqcongB}. Similarly, Theorem \ref{exp2} gives
	\begin{equation*}
		c_{p,\rho}^{(k)}(\alpha,a)
		=
		\sum_{m=0}^{p}
		\frac{(1)_{m,\rho}}{(\alpha m+a)^k}
		S_{1,\rho}(p,m).
	\end{equation*}
	The terms $m=0$ and $1\leq m\leq p-1$ vanish modulo $p$ since 
	\begin{equation*}
		S_{1,\rho}(p,m)\equiv 0\pmod p,
	\end{equation*}
	while the term $m=p$ is congruent to $\frac{(1)_{p,\rho}}{(\alpha p+a)^k}$ modulo $p$. Thus,
	\begin{equation*}
		c_{p,\rho}^{(k)}(\alpha,a)\equiv \frac{(1)_{p,\rho}}{(\alpha p+a)^k}\pmod p.
	\end{equation*}
	This proves \eqref{eqcongc}. Applying the same argument to
	Theorem \ref{exp3} gives
	\begin{equation*}
		\widehat{c}_{p,\rho}^{(k)}(\alpha,a)\equiv \frac{(-1)^{p}(1)_{p,\rho}}{(\alpha p+a)^k}\pmod p,
	\end{equation*}
	which proves \eqref{eqcongbc}.
\end{proof}

\subsection{Derivative Relations}
In this subsection, we examine the behavior of the degenerate two-parameter Hurwitz-Lerch type poly-Cauchy numbers of the first kind under differentiation of their generating function. By analyzing the derivative of the generating function expressed in terms of the generalized Hurwitz-Lerch-type polylogarithmic function, we derive a new sequence $D_{n,\rho}^{(k)}(\alpha,a)$ that encapsulates the rate of change of the poly-Cauchy sequence. The resulting identity involves Stirling numbers of the first kind and provides deeper insight into the analytic structure and differential properties of these generalized special numbers.

\begin{theorem}\label{thmderivativec}
	Let \(c_{n,\rho}^{(k)}(\alpha,a)\) be the degenerate two-parameter Hurwitz-Lerch type
	poly-Cauchy numbers of the first kind, defined by
	\begin{eqnarray*}
		\Phi_{f,\rho}(\log_{\rho}(1+t),k,\alpha,a)
		=
		\sum_{n=0}^{\infty}c_{n,\rho}^{(k)}(\alpha,a)\frac{t^n}{n!}.
	\end{eqnarray*}
	Then
	\begin{eqnarray*}
		\frac{d}{dt}\Phi_{f,\rho}(\log_{\rho}(1+t),k,\alpha,a)
		=
		(1+t)^{\rho-1}
		\sum_{n=0}^{\infty}D_{n,\rho}^{(k)}(\alpha,a)\frac{t^n}{n!},
	\end{eqnarray*}
	where
	\begin{eqnarray*}
		D_{n,\rho}^{(k)}(\alpha,a)
		=
		\sum_{m=1}^{n+1}
		\frac{(1)_{m,\rho}}{(\alpha m+a)^k}
		S_{1,\rho}(n,m-1).
	\end{eqnarray*}
\end{theorem}

\begin{proof}
	By (11), we have
	\begin{eqnarray*}
		\Phi_{f,\rho}(\log_{\rho}(1+t),k,\alpha,a)
		&=&
		\sum_{m=0}^{\infty}
		\frac{(1)_{m,\rho}(\log_{\rho}(1+t))^m}{m!(\alpha m+a)^k}.
	\end{eqnarray*}
	Differentiating both sides with respect to $t$, we obtain
	\begin{eqnarray*}
		\frac{d}{dt}\Phi_{f,\rho}(\log_{\rho}(1+t),k,\alpha,a)
		&=&
		\sum_{m=1}^{\infty}
		\frac{(1)_{m,\rho}}{m!(\alpha m+a)^k}
		m(\log_{\rho}(1+t))^{m-1}\frac{d}{dt}\log_{\rho}(1+t)\\
		&=&
		(1+t)^{\rho-1}
		\sum_{m=1}^{\infty}
		\frac{(1)_{m,\rho}(\log_{\rho}(1+t))^{m-1}}{(m-1)!(\alpha m+a)^k}.
	\end{eqnarray*}
	Let \(j=m-1\). Then
	\begin{eqnarray*}
		\frac{d}{dt}\Phi_{f,\rho}(\log_{\rho}(1+t),k,\alpha,a)
		&=&
		(1+t)^{\rho-1}
		\sum_{j=0}^{\infty}
		\frac{(1)_{j+1,\rho}(\log_{\rho}(1+t))^j}{j!(\alpha(j+1)+a)^k}.
	\end{eqnarray*}
	Using (13),
	\begin{eqnarray*}
		(\log_{\rho}(1+t))^j
		=
		j!\sum_{n=j}^{\infty}S_{1,\rho}(n,j)\frac{t^n}{n!}.
	\end{eqnarray*}
	Hence,
	\begin{eqnarray*}
		\frac{d}{dt}\Phi_{f,\rho}(\log_{\rho}(1+t),k,\alpha,a)
		&=&
		(1+t)^{\rho-1}
		\sum_{j=0}^{\infty}
		\frac{(1)_{j+1,\rho}}{(\alpha(j+1)+a)^k}
		\sum_{n=j}^{\infty}S_{1,\rho}(n,j)\frac{t^n}{n!}
		\\
		&=&
		(1+t)^{\rho-1}
		\sum_{n=0}^{\infty}
		\left(
		\sum_{j=0}^{n}
		\frac{(1)_{j+1,\rho}}{(\alpha(j+1)+a)^k}
		S_{1,\rho}(n,j)
		\right)\frac{t^n}{n!}.
	\end{eqnarray*}
	Replacing \(j+1\) by \(m\), we obtain
	\begin{eqnarray*}
		\frac{d}{dt}\Phi_{f,\rho}(\log_{\rho}(1+t),k,\alpha,a)
		&=&
		(1+t)^{\rho-1}
		\sum_{n=0}^{\infty}
		\left(
		\sum_{m=1}^{n+1}
		\frac{(1)_{m,\rho}}{(\alpha m+a)^k}
		S_{1,\rho}(n,m-1)
		\right)\frac{t^n}{n!}.
	\end{eqnarray*}
	This completes the proof.
\end{proof}

\begin{theorem}\label{thmderivativebc}
	Let $\widehat{c}_{n,\rho}^{(k)}(\alpha,a)$ be the degenerate two-parameter Hurwitz-Lerch type
	poly-Cauchy numbers of the second kind, defined by
	\begin{eqnarray*}
		\Phi_{f,\rho}(-\log_{\rho}(1+t),k,\alpha,a)
		=
		\sum_{n=0}^{\infty}bc_{n,\rho}^{(k)}(\alpha,a)\frac{t^n}{n!}.
	\end{eqnarray*}
	Then
	\begin{eqnarray*}
		\frac{d}{dt}\Phi_{f,\rho}(-\log_{\rho}(1+t),k,\alpha,a)
		=
		(1+t)^{\rho-1}
		\sum_{n=0}^{\infty}\widehat{D}_{n,\rho}^{(k)}(\alpha,a)\frac{t^n}{n!},
	\end{eqnarray*}
	where
	\begin{eqnarray*}
		\widehat{D}_{n,\rho}^{(k)}(\alpha,a)
		=
		\sum_{m=1}^{n+1}
		(-1)^m
		\frac{(1)_{m,\rho}}{(\alpha m+a)^k}
		S_{1,\rho}(n,m-1).
	\end{eqnarray*}
\end{theorem}

\begin{proof}
Analogous to Theorem \ref{thmderivativec}.
\end{proof}

\begin{theorem}\label{thmderivativeB}
	Let \(B_{n,\rho}^{(k)}(\alpha,a)\) be the degenerate two-parameter Hurwitz-Lerch type
	poly-Bernoulli numbers, defined by
	\begin{eqnarray*}
		\Phi_{\rho}(1-e_{\rho}(-t),k,\alpha,a)
		=
		\sum_{n=0}^{\infty}B_{n,\rho}^{(k)}(\alpha,a)\frac{t^n}{n!}.
	\end{eqnarray*}
	Then
	\begin{eqnarray*}
		\frac{d}{dt}\Phi_{\rho}(1-e_{\rho}(-t),k,\alpha,a)
		=
		(1-\rho t)^{\frac{1}{\rho}-1}
		\sum_{n=0}^{\infty}E_{n,\rho}^{(k)}(\alpha,a)\frac{t^n}{n!},
	\end{eqnarray*}
	where
	\begin{eqnarray*}
		E_{n,\rho}^{(k)}(\alpha,a)
		=
		(-1)^n
		\sum_{m=1}^{n+1}
		(-1)^{m-1}
		\frac{(1)_{m,\rho}m!}{(\alpha m+a)^k}
		S_{2,\rho}(n,m-1).
	\end{eqnarray*}
\end{theorem}

\begin{proof}
Analogous to Theorem \ref{thmderivativec}	.
\end{proof}
\subsection{Conclusion}
In this paper, we introduced degenerate Hurwitz-Lerch type poly-Bernoulli and poly-Cauchy
numbers with parameters $(\alpha,a)$ by means of degenerate Hurwitz-Lerch type zeta
functions and their factorial analogues. We established generating functions for these families and derived explicit formulas in terms of the degenerate Stirling numbers of the
first and second kinds.

We further obtained orthogonality relations together with duality and inversion identities.  These
clarify the structural connections among the degenerate poly-Bernoulli numbers and the
corresponding poly-Cauchy numbers of the first and second kinds. Moreover, we derived
congruence properties modulo a prime and established derivative relations for the associated
generating functions, thereby revealing both arithmetic and analytic aspects of these
generalized sequences.

Taken together, these results extend several known properties of classical and generalized
poly-Bernoulli and poly-Cauchy numbers to the degenerate two-parameter Hurwitz-Lerch
setting. They also provide a useful framework for further investigations in combinatorics,
special functions, and analytic number theory.

Possible directions for future work include the study of recurrence relations, asymptotic
behavior, and q-analogues of these families as well as further connections with generalized
Stirling numbers and other Hurwitz-Lerch type special functions.


\begin{thebibliography}{99}
	
	\bibitem{apostol}
	T. M. Apostol, On the Lerch zeta function, \textit{Pacific J. Math.} \textbf{1}(2)
	(1951), 161-167.
	
	\bibitem{arakawa}
	T. Arakawa and M. Kaneko, Multiple zeta values, poly-Bernoulli numbers, and
	related zeta functions, \textit{Nagoya Math. J.} \textbf{153} (1999), 189-209.
	
	\bibitem{carlitz}
	L. Carlitz, Degenerate Stirling, Bernoulli and Eulerian numbers,
	\textit{Utilitas Math.} \textbf{15} (1979), 51-88.
	
	\bibitem{cenkci}
	M. Cenkci and P. T. Young, Generalizations of poly-Bernoulli and poly-Cauchy
	numbers, \textit{Eur. J. Math.} \textbf{1} (2015), 799-828.
	doi:10.1007/s40879-015-0071-3.
	
	\bibitem{comtet}
	L. Comtet, \textit{Advanced Combinatorics}, D. Reidel Publishing Co.,
	Dordrecht, 1974.
	
	\bibitem{corcino1}
	C. B. Corcino, R. B. Corcino, J. M. Ontolan, and C. M. Perez-Fernandez,
	The Hankel transform of $q$-noncentral Bell numbers,
	\textit{Int. J. Math. Math. Sci.} \textbf{2015} (2015), Article ID 417327,
	10 pages. doi:10.1155/2015/417327.
	
	\bibitem{corcino2}
	R. B. Corcino and C. B. Corcino, The Hankel transform of generalized Bell
	numbers and its $q$-analogue, \textit{Util. Math.} \textbf{89} (2012), 297-309.
	
	\bibitem{corcino3}
	R. B. Corcino and C. B. Montero, A $q$-analogue of Rucinski-Voigt numbers,
	\textit{ISRN Discrete Math.} \textbf{2012} (2012), Article ID 592818, 18 pages.
	doi:10.5402/2012/592818.
	
	\bibitem{corcino4}
	R. B. Corcino and C. Barrientos, Some theorems on the $q$-analogue of the
	generalized Stirling numbers, \textit{Bull. Malays. Math. Sci. Soc.}
	\textbf{34}(3) (2011), 487-501.
	
	\bibitem{corcino5}
	R. B. Corcino, L. C. Hsu, and E. L. Tan, A $q$-analogue of generalized
	Stirling numbers, \textit{Fibonacci Quart.} \textbf{44}(2) (2006), 154-165.
	
	\bibitem{corcino-herrera}
	R. B. Corcino and M. L. Herrera, Some properties of the limit of the differences
	of the generalized factorial, \textit{Matimyas Matematika} \textbf{32}(2)
	(2009), 21-28.
	
	\bibitem{corcino-limit}
	R. B. Corcino, C. B. Corcino, and M. L. Herrera, On the limiting form of the
	differences of the generalized factorial, preprint.
	
	\bibitem{graham}
	R. L. Graham, D. E. Knuth, and O. Patashnik, \textit{Concrete Mathematics:
		A Foundation for Computer Science}, 2nd ed., Addison-Wesley, Reading, MA, 1994.
	
	\bibitem{kaneko}
	M. Kaneko, Poly-Bernoulli numbers, \textit{J. Th\'eor. Nombres Bordeaux}
	\textbf{9}(1) (1997), 221-228.
	
	\bibitem{Kim-qzeta}
	T. Kim, A new approach to $q$-zeta function, arXiv:math/0502005 [math.NT],
	2005. doi:10.48550/arXiv.math/0502005.
	
	\bibitem{Kim-lambda}
	T. Kim, $\lambda$-analogue of Stirling numbers of the first kind,
	\textit{Adv. Stud. Contemp. Math. (Kyungshang)} \textbf{27}(3) (2017),
	423-429.
	
	\bibitem{Kim-2}
	T. Kim, A note on degenerate Stirling polynomials of the second kind,
	\textit{Proc. Jangjeon Math. Soc.} \textbf{20}(3) (2017), 319-331.
	
	\bibitem{KimBell}
	T. Kim, Degenerate complete Bell polynomials and numbers,
	\textit{Proc. Jangjeon Math. Soc.} \textbf{20}(4) (2017), 533-543.
	
	\bibitem{Kim-1-1}
	D. S. Kim and T. Kim, A note on a new type of degenerate Bernoulli numbers,
	\textit{Russian J. Math. Phys.} \textbf{27}(2) (2020), 227-235.
	doi:10.1134/S1061920820020090.
	
	\bibitem{KimKim2018}
	T. Kim and D. S. Kim, A note on degenerate Stirling numbers of the first kind,
	arXiv:1802.00896 [math.NT], 2018.
	
	\bibitem{KimKimGamma}
	T. Kim and D. S. Kim, Degenerate Laplace transform and degenerate gamma function,
	\textit{Russian J. Math. Phys.} \textbf{24}(2) (2017), 241-248.
	
	\bibitem{KimKimType2}
	T. Kim and D. S. Kim, A note on Type 2 Changhee and Daehee polynomials,
	\textit{Rev. R. Acad. Cienc. Exactas F\'is. Nat. Ser. A Mat. RACSAM}
	\textbf{113} (2019), 2783-2791.
	
	\bibitem{KimKimSymmetry}
	T. Kim and D. S. Kim, Some relations of two Type 2 polynomials and discrete
	harmonic numbers and polynomials, \textit{Symmetry} \textbf{12}(6) (2020),
	Article 905. doi:10.3390/sym12060905.
	
	\bibitem{KimKimJang}
	T. Kim, D. S. Kim, and G.-W. Jang, Extended Stirling polynomials of the second
	kind and extended Bell polynomials, \textit{Proc. Jangjeon Math. Soc.}
	\textbf{20}(3) (2017), 365-376.
	
	\bibitem{KimDolgy}
	T. Kim, D. S. Kim, and D. V. Dolgy, On partially degenerate Bell numbers and
	polynomials, \textit{Proc. Jangjeon Math. Soc.} \textbf{20}(3) (2017), 337-345.
	
	\bibitem{Kim-3}
	T. Kim, D. S. Kim, H. Y. Kim, and J. Kwon, Degenerate Stirling polynomials
	of the second kind and some applications, \textit{Symmetry} \textbf{11}(8)
	(2019), Article 1046. doi:10.3390/sym11081046.
	
	\bibitem{KimPolyBernoulli}
	T. Kim, D. S. Kim, H. Y. Kim, and L.-C. Jang, Degenerate poly-Bernoulli
	numbers and polynomials, \textit{Informatica} \textbf{31}(3) (2020), 1-7.
	
	\bibitem{KimJindalrae}
	T. Kim, D. S. Kim, L.-C. Jang, and H. Lee, Jindalrae and Gaenari numbers and
	polynomials in connection with Jindalrae-Stirling numbers,
	\textit{Adv. Difference Equ.} \textbf{2020} (2020), Article No. 245, 19 pages.
	doi:10.1186/s13662-020-02701-1.
	
	\bibitem{komatsu}
	T. Komatsu, Poly-Cauchy numbers, \textit{Kyushu J. Math.} \textbf{67}(1)
	(2013), 143-153.
	
	\bibitem{lacpao1}
	N. B. Lacpao, R. Corcino, and M. A. R. Vega, Hurwitz-Lerch type
	multi-poly-Cauchy numbers, \textit{Mathematics} \textbf{7}(4) (2019),
	Article 335. doi:10.3390/math7040335.
	
	\bibitem{lacpao2}
	N. B. Lacpao, Some explicit formulas of Hurwitz-Lerch type poly-Cauchy
	polynomials and poly-Bernoulli polynomials, \textit{European J. Pure Appl. Math.}
	\textbf{16}(3) (2023), 1747-1761. doi:10.29020/nybg.ejpam.v16i3.4825.
	
	\bibitem{lewin}
	L. Lewin, \textit{Polylogarithms and Associated Functions}, North-Holland,
	New York, 1981.
	
	\bibitem{ParkDaehee}
	J.-W. Park, B. M. Kim, and J. Kwon, On a modified degenerate Daehee polynomials
	and numbers, \textit{J. Nonlinear Sci. Appl.} \textbf{10}(3) (2017),
	1108-1115.
	
	\bibitem{roman}
	S. Roman, \textit{The Umbral Calculus}, Pure and Applied Mathematics, Vol. 111,
	Academic Press, New York, 1984.
	
	\bibitem{simsek}
	Y. Simsek, Identities on the Changhee numbers and Apostol-type Daehee
	polynomials, \textit{Adv. Stud. Contemp. Math. (Kyungshang)} \textbf{27}(2)
	(2017), 199-212.
	
	\bibitem{srivastava}
	H. M. Srivastava and J. Choi, \textit{Series Associated with the Zeta and Related
		Functions}, Kluwer Academic Publishers, Dordrecht, 2001.
	
	\bibitem{srivastava-1}
	H. M. Srivastava and J. Choi, \textit{Zeta and $q$-Zeta Functions and Associated
		Series and Integrals}, Elsevier, Amsterdam, 2012.
	
	\bibitem{YuStirling}
	H. Yu, A generalization of Stirling numbers, \textit{The Fibonacci Quarterly}
	\textbf{36}(3) (1998), 252--258.
	
\end{thebibliography}
\end{document}